\documentclass[10pt]{article}
\usepackage{a4wide}
\usepackage[T1]{fontenc}
\usepackage{multicol}% pour ecrire en double colonne
\usepackage{latexsym}
\usepackage{amssymb}
\usepackage{amsmath}

\newcommand\aaa{\end{colonnes}\end{document}}
\newcommand{\diag}{\rm  \diag}
\newcommand{\R}{\mathbb R}

\DeclareRobustCommand{\qed}{%
\ifmmode % if math mode, assume display: omit penalty etc.
\else \leavevmode\unskip\penalty9999 \hbox{}\nobreak\hfill
\fi
\quad\hbox{\openbox}}
\newcommand\un[1]{\,\rlap{{1}}\kern.22em \mbox{l}_{#1}}

\newcommand{\openbox}{\leavevmode
  \hbox to.77778em{%
  \hfil\vrule
  \vbox to.675em{\hrule width.6em\vfil\hrule}%
  \vrule\hfil}}
\def\1{{\rm 1\mskip-4,4mu l}}

\newenvironment{colonnes}{%
       \columnsep 8mm\begin{multicols}{2}}{%
       \end{multicols}
                          }

\newcommand\ad[1]{\mbox{\rm ad}_{#1}}

\begin{document}
\pagestyle{empty}

\rule{1cm}{0cm}
\vskip 4.4cm
\begin{center}
{\bf GEODESICS IN NILPOTENT LIE GROUPS}\\
\rule{1cm}{0cm}\\
{\bf Renée Abib}
\end{center}
\begin{center}
{\em Laboratoire de Mathématiques Raphaël Salem - UMR6085-CNRS \\
Université de Rouen \\
76821 Mont Saint Aignan Cedex, France} \\
{\em E-mail: Renee.Abib@univ-rouen.fr}
\end{center}
\vspace{0,5cm}

\begin{center}\begin{minipage}{125mm}
\noindent Abstract: We study the sub-Riemannian geodesic problem in
the $2n+1$-Heisenberg group. An independent set of $2n$ left-invariant
vector fields determines a $SR$ structure. We are interested in
sub-Riemannian length minimizers for this structure. We use the
hamiltonian formalism and apply Pontryagin maximum principle to write
the equations for the extremals and describe some properties of the
extremal curves. We consider geodesics in the group with a
left-invariant Riemannian metric. We obtain their equations and prove
that the set of directions of all rays issuing form zero is the sphere
of unit horizontal directions. \\ 
\rule{1cm}{0cm}
\end{minipage}\end{center}
\rule{1cm}{0cm}\\
\begin{colonnes}
\begin{center}
1. INTRODUCTION
\end{center}
Let $L$ be the $2n+1$--dimensional nilpotent Lie algebra with the
following commutation rules in some basis
$X_1,Y_1,X_2,Y_2,\dots,X_n,Y_n,T$: 
\begin{align*}
&[X_i,Y_j]=2\delta_{ij}T, \\
&[X_i,X_j]=[Y_i,Y_j]=\ad{T}=0,
\end{align*}
and $G$ be the connected, simply connected Lie group with the Lie
algebra $L$. We consider $X_i,Y_i,T$ as left-invariant vector fields
on $G$. The set $\{X_1,Y_1,\cdots,X_n,Y_n\}$ determines a
left-invariant sub-Riemannian structure on $G$: \\  

\indent $D=\hbox{span }(X_1,Y_1,\cdots,X_n,Y_n),  \left\langle X_i,Y_j\right\rangle$ \\ 
\indent $=\left\langle  X_i,X_j\right\rangle =\left\langle
  Y_i,Y_j\right\rangle=\delta_{ij}              \, i,j=1,2,\cdots,n,$
\\ 

\noindent called the flat $(2n,2n+1)$ sub-Riemannian structure. Such
structure is unique, up to isomorphism of Lie groups. We are
interested in sub-Riemannian length minimizers for this sub-Riemannian
structure, i.e., in solutions to the following optimal control
problem: 
\begin{equation} \label{equa1} 
\dot{\gamma}=\sum_{i=1}^n \left(u_i X_i(\gamma)+v_i Y_i(\gamma)\right)
\end{equation}
\hspace{4,5cm} $\gamma \in G, u_i,v_i \in \R,$ \\
\begin{equation} \label{equa2}
\gamma(0)=q_0, \gamma(T)=q_1 \text{ fixed},
\end{equation}
\begin{equation} \label{equa3}
l=\int_0^T\sqrt{\sum_{i=1}^n (u_i^2+v^2_i)}\rightarrow \min.
\end{equation}
Problem (\ref{equa1})-(\ref{equa3}) was considered by Brockett(1981),
and Liu and Sussmann (1995) for $n=1$. $G$ is a $2$-step nilpotent Lie
group and $L=D\oplus[L,L]$. In this note we consider also geodesic
lines in the Lie group $G$ with a left-invariant Riemannian metric. \\ 
\\
\noindent The flat $(2n,2n+1)$ sub-Riemannian structure gives a local
nilpotent approximation for an arbitrary sub-Riemannian structure with
growth vector $(2n,2n+1)$, see Montgomery (2002). The dynamics of a
classical electric charge in the plan under the influence of a
perpendicular magnetic field can be described by means of flat
sub-Riemannian structure. \\ 

\noindent In general, if $G$ is a $2$-step nilpotent Lie group and
$(D,\left\langle ,\right\rangle)$ is a left-invariant sub-Riemannian
structure such that $L=D\oplus[L,L]$, any minimizer is normal in some
sub-group of $G$ and hence any minimizer is smooth, see Abib
(2002). This extends known results on the so-called Gaveau-Brockett
problem, Brockett (1981) and Gaveau (1977). 

\begin{center}
2. MODEL
\end{center}
We choose the following model for the flat $(2n,2n+1)$ sub-Riemannian structure: \\
\\
$G=\R^{2n+1}$ with coordinates $q=(x,z)$ where
$x=(x_1,y_1,\dots\dots,x_n,y_n),z\in \R$, and group law defined as
follows: 
\begin{align*}
&(x,z)(x',z')
=(x+x',z+z'+\rule{2cm}{0cm}\\
&\phantom{(x,z)(x',z')=}+\sum_{i=1}^nx_iy'_i-\sum_{i=1}^nx'_iy_i)\\
&X_i
=\frac{\partial}{\partial
    x_i}-y_i\frac{\partial}{\partial z}\\
&Y_i=\frac{\partial}{\partial y_i}+x_i\frac{\partial}{\partial z}\\
&D=\hbox{Ker } \theta\\
&\theta=dz+\sum_{i=1}^n(y_idx_i-x_idy_i)
\end{align*}

Control system (\ref{equa1}) has the full rank and the state space $G$
is connected, thus the systems is globally controllable on $G$. A
standard existence theorem from optimal control theory, see Liu and
Sussmann (1995), implies that any point can be connected with $q_0 \in
G$ by an admissible length minimizer. \\ 
\\
\begin{center}
3. EXTREMALS
\end{center}
A admissible curve for $D$ is an absolutely continuous curve on $G$ which is tangent to $D$ almost everywhere; applying the Pontryagin maximum principle each minimizer parametrized by arc-length is a normal extremal or abnormal extremal, see Liu and Sussmann (1995).\\
\\
\noindent 3.1 {\em Abnormal extremals}\\
\\
An abnormal extremal is a admissible curve for $D$ which is the projection onto $G$ of characteristic curve (in the symplectic sense) of the annihilator $D^{\bot}$ of $D$ in $T^* G$. Introduce the linear Hamiltonians corresponding to the basis fields: 
$$h_i(\lambda)=\lambda(X_i),\,  k_i(\lambda)=\lambda(Y_i),\,  h(\lambda)=\lambda(T),$$  
$$\lambda \in T^* G.$$
$D^{\bot}$ is defined by equations $h_i=k_i=0$, $i=1,\cdots,n$. Assume $\Gamma(t) \in T^* G-\{0\}$ caracteristic curve of $D^{\bot}$, then $\dot{\Gamma}(t) \in \hbox{span } \{\vec{h_i}(\Gamma(t)), \, \vec{k_i} (\Gamma(t)$ ; $i=1,2,\cdots,n\}$ and $h_i(\Gamma(t))=k_i(\Gamma(t))=0$.\\
Abnormal extremals are exactly the constant curves.\\
\\
\noindent 3.2 {\em Normal extremals} \\
\\
Normal biextremals are trajectories of the Hamiltonian system
\begin{equation} \label{equa4}
\dot{\lambda}=\vec{H}(\lambda), \, \lambda \in T^* G
\end{equation}
with the sub-Riemannian Hamiltonian $H=\frac{1}{2}\sum_{i=1}^n(k_i^2+k_i^2)$. In the coordinates $(q,\lambda)=(x_1,\cdots,x_n,y_1,\cdots,y_n,z,\xi_1,\cdots,\xi_n,\eta_1,\cdots,\eta_n,\zeta)$ of $T^*G$ system (\ref{equa4}) reads 
\begin{equation} \label{equa5}
\dot{x}_i=\xi_i+y_i \zeta
\end{equation}
\begin{equation} \label{equa6}
\dot{y}_i=\eta_i-x_i \zeta
\end{equation}
\begin{equation} \label{equa7}
\dot{z}=\sum_{i=1}^n(\dot{x}_i y_i-\dot{y}_i x_i)
\end{equation}
\begin{equation} \label{equa8}
\dot{\xi_i}=\dot{y_i} \zeta
\end{equation}
\begin{equation} \label{equa9}
\dot{\eta}_i=\dot{x_i} \zeta
\end{equation}
\begin{equation} \label{equa10}
\dot{\zeta}=0
\end{equation}
If the constant $\zeta$ is $0$, the curve $(x(t),y(t))$ is the
straight line, along which the motion takes place with velocity
constant. The function $z(t)$ is then computed by solving
(\ref{equa7}) and the curve $(x(t), y(t), z(t))$ is a straight line in
$G$. If $\zeta$ is non zero constant, the equations
(\ref{equa8})-(\ref{equa9}) imply 
\begin{align*}
&\ddot{x_i}-2\dot{y_i} \zeta=0,
&\ddot{y_i}+2\dot{x_i} \zeta=0. 
\end{align*}
So 
\begin{align*}
&\dot{x_i}=A_i \cos 2\zeta t+B_i \sin 2\zeta t, \, \\
&\dot{y_i}=-A_i \sin 2 \zeta t+B_i \cos 2 \zeta t,\\
&A_i,B_i
\end{align*} 
are constant. The condition that our curve
$q(t)=(x(t),y(t),z(t))$ is parametrized by arc-length then says that
$\sum_{i=1}^n (A_i^2+B_i^2)=1$. We obtain
\begin{align*}
&x_i(t)=\frac{r_i}{2\zeta}(\cos \theta_i-\cos (2\zeta t+\theta_i)),\\
&y_i(t)=\frac{r_i}{2\zeta}(\sin(2\zeta t+\theta_i)-\sin \theta_i)
\end{align*} 
with $A_i=r_i \sin \theta_i, B_i=r_i \cos \theta_i$. The equation
(\ref{equa7}) imply 
\begin{align*}
&\ddot{z}=\zeta \sum_{i=1}^n \frac{d}{d t}(x_i^2+y_i^2) \\
\text{and then }&z(t)=\sum_{i=1}^n r_i^2 (\frac{t}{2
  \zeta}-\frac{\sin 2 \zeta t}{(2 \zeta)^2}). 
\end{align*} 
The Pontryagin maximum
principle, applied to this case,  implies that all the minimizers
curve parametrized by arc-length are normal extremal $q(t)$. 
\\

\begin{center} 
4. CONTACT TRANSFORMATIONS
\end{center}
For $n=1$, the Lie algebra of symmetries of the flat (2,3)
distribution $D$ is parametrized by arbitrary smooth functions of
three variables and the Lie algebra of symmetries of the flat (2,3)
sub-Riemannian structure is $4$-dimensional, see Sachkov (1998). 
\\
\\
\noindent For $n \geq 1$, let $C$ be the sheaf of germs of infinitesimal
automorphisms $X$ such that $L_X \theta=f \theta$, the function $f$
depending on $X$. $C$ is the sheaf of contact transformations.  
\\
\\
Let $X \in C$ vanishing at $0$, and $\phi_t$ be the one parameter
group generated by $X$. Then $\phi_t(0)=0$ and $\phi_t^* \theta=f_t
\theta$. Let $V_1$ be the subspace of $L$ defined by
$\theta_0=0$. Then $\phi_t$ leaves $V_1$ invariant. Furthermore
$\phi_t^* d \theta=df_t \wedge \theta + f_t d\theta$. This implies
that the form $d \theta$ restricted to $V_1$ is preserved up to scalar
factor. If we pass to the infinitesimal we see that the linear
isotropy algebra of $C$ is a subalgebra of the algebra $g$ defined as
follows:  
\\
\\
let $u \in L$ not lying in $V_1$ and $g=K+N+M$ 
where $K$ is the set of all $A \in \hbox{End}(L)$ with $\hbox{A}u=0$ and $A/V_1 \in \hbox{sp}(V_1)$, \\
$N$ is the of all $B$ of the form $\hbox{B}v_1=0$ for $v_1 \in V_1$ and $\hbox{B}u\in V_1$, \\
$M$ is the set of all multiples of the linear transformation $C$ where $C v_1=v_1$ for $v_1 \in V_1$ and $Cu=2u$.\\
 
\noindent The space $g^{(1)} \subset \hbox{Hom } (g,L)$ called the
first prolongation of $g$, is the set of all $T \in \hbox{Hom }(g,L)$
which satisfy $(T(\omega))(v)=(T(v))(\omega)$ for all $\omega$, $v \in
L$. In the case, $g^{(1)}=\hbox{sp }(V_1)^{(1)}+g$ and
$g^{(k)}=\hbox{sp  }(V_1)^{(k)} + g^{(k-1)}$ for all $k$. The algebra
$\hbox{sp }(V_1)$ is of infinite type because $\hbox{sp }(V_1)^{(k)}$
can be identified with $S^{k+2}(V_1^*)$; then $g$ is of infinite type.
\\ 

\noindent Corresponding to the coordinates $z,x_1,\cdots,x_n, y_1$, 
$\cdots, y_n$, let us choose the basis
$u=e_0,e_1,\cdots,e_n,e_{n+1},\cdots,e_2$ of $L$. Then it is easy to
verify that the following vector fields are all infinitesimal contact
transformations (where $\theta=dz+\sum_{i=1}^n(y_i dx_i - x_i dy_i)$):
\\ 
\\
$\alpha)$ ${\displaystyle \frac{\partial}{\partial x_i} + y_i
  \frac{\partial}{\partial z}, \, \frac{\partial}{\partial y_i}-x_i
  \frac{\partial}{\partial z}, \frac{\partial}{\partial z}}$ \\ 
\\
$\beta)$ ${\displaystyle \sum_{i,j=1}^n \Bigg[A_{ij} \, x_i
  \frac{\partial}{\partial x_j} + A_{i+n,j} y_i
  \frac{\partial}{\partial x_j}+A_{i,j+n} x_i \frac{\partial}{\partial
    y_j}}$ \\ 

${\displaystyle +A_{i+n,j+n} y_i \frac{\partial}{\partial y_j}\Bigg]}$ \\

where ${\displaystyle \sum_{r,s=1}^{2n} A_{rs} \, e_s \otimes e^*_r \, \in sp(V_1)},$ \\
\\
${\displaystyle 2z \frac{\partial}{\partial z}+\sum_{j=1}^n x_j \frac{\partial}{\partial x_j} +\sum_{j=1}^n \, y_j \frac{\partial}{\partial y_j}},$ \\ 

${\displaystyle z \frac{\partial}{\partial x_i}+y_i \Bigg[\sum_{j=1}^n \Bigg(x_j\frac{\partial}{\partial x_j}+y_j\frac{\partial}{\partial y_j}\Bigg)+z\frac{\partial}{\partial z}\Bigg],}$\\
\\
${\displaystyle z\frac{\partial}{\partial y_i}+x_i\Bigg[\sum_{j=1}^n \Bigg(x_j\frac{\partial}{\partial x_j}+y_j\frac{\partial}{\partial y_j}\Bigg)+z\frac{\partial}{\partial z}\Bigg]}$ \\
\\
$\gamma)$ ${\displaystyle z\Bigg[\sum_{j=1}^n\Bigg(x_j\frac{\partial}{\partial x_j}+y_j\frac{\partial}{\partial y_j}\Bigg)+z\frac{\partial}{\partial z}\Bigg]}$ \\
\\
\\
>From $\alpha)$ we see that $C$ is transitive. From $\beta)$ we see that the linear isotropy algebra of $C$ is indeed $g$. A direct computation shows that the vector fields $\alpha)$, $\beta)$ and $\gamma)$ form a Lie algebra. $C$ is an {\em non-flat, infinite, transitive} LAS (Lie algebra sheaf).
\\
\begin{center}
5. RIEMANNIAN CASE
\end{center}
We define the left-invariant metric on $G$ by taking $X_i$, $Y_i$, $T$
as the orthonormal frame. Let $\nabla$ the Riemannian connection of
left-invariant metric an
$\{e_h,h=1,2,\cdots,2n+1\}=\{X_1,Y_1,\cdots,X_n,Y_n,T\}$. Then
$\nabla_{e_i}e_j=\frac{1}{2}\sum_{k=1}^{2n}(c_{ij}^k+c_{ki}^j-c_{jk}^i)e_k$
where $[e_i,e_j]=\sum_{k=1}^{2n}c_{ij}^k e_k$.\\ 
\\
 For $G$ the matrice $(\nabla_{e_i}^{e_j})_{i,j}$ is: 
\begin{align} \label{equa11}
\left(
\begin{array}{cccccc} 
0 & T & 0 & \cdots & 0 & -Y_1 \\  
-T & 0 & 0 &  \cdots & 0 & X_1 \\
\vdots &  &   &  & \vdots & \vdots \\
0 & \cdots & \cdots & 0 & T & -Y_n \\
0 & \cdots & 0 & -T & 0 & X_n \\
-Y_1 & X_1 & \cdots & -Y_n & X_n & 0
\end{array}
\right)
\end{align}
Let
$\dot{c}(t)=\sum_{i=1}^n\left(u_i(t)X_i(t)+v_i(t)Y_i(t)\right)+\gamma(t)T$
geodesics issuing from $0 \in G$. Then
$\nabla_{\dot{c}(t)}^{\dot{c}(t)}=0$ and the table (\ref{equa11}) give
\begin{align*}
&\dot{u}_i+2\gamma v_i=0,\\ 
&\dot{v}_i-2\gamma u_i=0,\\
&\dot{\gamma}=0. 
\end{align*}
Therefore, because the parameter $t$ is natural we
have $\sum_{i=1}^n(u_i^2+v_i^2)+\gamma^2=1$ and we could take
$\gamma(t)=\gamma$ where 
the constant $\gamma\in[-1,1]$  
%$|\gamma| \leq 1$ 
is the cosine of the angle
between $\dot{c}(0)$ and the $T$-axe. \\ 
\\
 For $\gamma \neq 0$, $u_i(t)=r_i\cos (2\gamma t+\theta_i)$ and
 $v_i(t)=r_i \cos (2 \gamma t+\theta_i)$ where $r_i^2=u^2_i+v_i^2$. \\
\\
In
 coordinates $(x_1,y_1,\cdots,x_n,y_n,z)$ the equations for geodesics
 $c(t)=(x_1(t),\cdots,y_n(t),z(t))$ are: 
\begin{equation*}
\dot{x}_i=u_i, \, \dot{y}_i=v_i,
\end{equation*}
\begin{equation} \label{equa12}
\dot{z}=\gamma+\sum_{i=1}^n v_i(t)x_i(t)-\sum_{i=1}^n u_i(t)y_i(t)
\end{equation}
We have \\

${\displaystyle x_i(t) = -\frac{r_i}{2\gamma}\cos 2\gamma t + a_i}$ \\

${\displaystyle y_i(t) = \frac{r_i}{2\gamma} \sin 2\gamma t + b_i}$ \\

${\displaystyle z(t) = \frac{3\gamma^2-1}{2\gamma}t+\frac{b_i}{2\gamma}r_i \cos 2\gamma t}$ \\

$\hspace{4cm} {\displaystyle +\frac{a_i}{2\gamma}r_i \sin 2\gamma t +c}$ \\

\noindent for some numbers $a_i$, $b_i$, $c$ which could be defined from the initial condition $c(0)=0$. \\
If $\gamma=0$, then they are horizontal and satisfy to the following:
$$
x_i(t)=u_i(0)t, \, y_i(t)=v_i(0)t, \, z(t)=0
$$

\noindent The non-horizontal geodesics (i.e., with $\gamma \neq 0$)
are not rays. Indeed. Let $|\gamma| \neq 1.$ For every geodesic $c(t)$
their $x_i$, $y_i$ coordinates are periodic functions. Their
projection on $z=0$ hyperspace is a circle. Take two points on $c(t)$
with difference of their parameters equals $\frac{\pi}{\gamma}$. Then,
besides the geodesics $c(t)$ their is also the vertical geodesic
$z(s)=(0,\cdots,0,s)$ connecting them. Its length is equal to the
difference of $z$ coordinates of considering points. But, as easy to
check this is always strictly less than $\frac{\pi}{\gamma}$, i.e.,
the length of the interval of the $c(t)$ between them, because $t$ is
natural one by definition. Therefore, every geodesic $c(t)$ with
$|\gamma|\neq 1$ is not minimal. Vertical geodesic $z(s)$ (with
$|\gamma|=1$) also is not minimal, because that the vector field
$Z(s)=X_1 \cos s + Y_1 \sin$ is parallel along it and the sectional
curvatures of $G$ in $2$-dimensional direction generated by $Z(s)$ and
$\dot{z}(s)$ equals $1$. This means that the index of every interval
of $z(s)$ with a length greater than $\pi$ is positive, and $z$ is not
minimal. \\ 

\noindent Now let $u s$ check that all horizontal geodesics (with
$\gamma=0$) are rays. Suppose that the horizontal geodesic $a(t)$
intersects in some point with another geodesic $c(s)$ issuing the same
point $0$. Because horizontal geodesics do not intersect, $c(s)$ is
non-horizontal one. For $c(t)$ from (\ref{equa12}) we see that the
length $\tilde{l}(c)$ of the projection curve
$\tilde{c}(t)=(c_1(t),\cdots,c_{2n}(t),0)$, measured in euclidean
coordinates, is strictly less than $t$, i.e., the length of the curve
$c$ in the left-invariant metric. Than in the submanifold $z=0$
considered as euclidean space the length of every chord, connecting
some points of $\tilde{c}$ is less than $\tilde{l}(c)$. If the
geodesic $c$ issuing from $0$ intersects some horizontal geodesic
$a(t)$ also issuing from $0$, then a is the chord of the corresponding
$\tilde{c}$. For horizontal geodesics their lengths in the
left-invariant metric coincide with the usual euclidean length in
$z=0$. The argument above means that the length of $a$ is strictly
less then that of $c$. Therefore, we have the following statement: the
set of directions of all rays issuing from $0$ in the group $G$ is the
sphere of unit horizontal directions in the point $0$ i.e.,
$\{w=(w_1,\cdots,w_{2n},0);\, \|w\|=1\}$.  
\vspace{0,4cm}
\begin{center}
REFERENCES
\end{center}

\noindent Abib, O.R. (2002). Sub-Riemannian geodesics \break 
\indent on Lie groups. In {\em Proceedings of the I Collo-\break 
\indent quium on Lie theory and application} (I. Bajo\break 
\indent and  E. SanMartin, Ed.), pp.3-12 , Vigo,\break 
\indent Spain.

\noindent Brockett, R. (1981). Control theory and singular\break 
\indent  Riemannian geometry. In: {\em New Directions in\break 
\indent Applied Mathematics} (P. Hilton and G. \break 
\indent Young, Ed.), pp.11-27. Springer-Verlag, New\break 
\indent York.

\noindent Gaveau, B. (1977). Principe de moindre action,\break 
\indent propagation de la chaleur et estimées sous\break
\indent elliptiques  sur certains groupes nilpotents.\break 
\indent {\em Acta Math.} 139, pp.94-153.

\noindent Liu, W. and H. Sussmann (1995). Shortest paths 

for sub-Riemannian metrics on rank-two dis-\break
\indent tributions. {\em Mem. Amer. Math. Soc.} 118,\break
\indent No 564.

\noindent Montgomery, R. (2002). A tour of Sub-Rieman-\break
\indent nian geometries, their geodesics and applica-\break
\indent tions. {\em American Mathematical Society}, \break
\indent Rhode Island.

\noindent Sachkov, Y. (1998). Symmetries of flat rank two\break 
\indent distributions and sub-Riemannian struc-\break  
\indent tures. {\em Technical report}, 98-151, Laboratoire\break 
\indent de Topologie de Dijon, France.

\end{colonnes}
\end{document}